\documentclass[11pt,twoside]{article}
\usepackage{mathrsfs}
\usepackage{graphics}
\usepackage{amssymb,amsmath}
\usepackage{graphicx}
\usepackage{amsmath}
\usepackage{amsthm}
\usepackage{amsfonts}
\usepackage{amssymb}
\usepackage{latexsym}
\usepackage{amscd,amsxtra}
\usepackage{subfigure}
\usepackage{amscd}
\usepackage{slashed}
\usepackage{titlesec}          

\renewcommand\thesection{\arabic{section}}
\titleformat{\section}{\normalfont\large\bfseries}{\thesection.}{0.5em}{}

\makeatletter  \evensidemargin0.4in \oddsidemargin0.4in

\newtheorem{theorem}{Theorem}[section]
\newtheorem{proposition}{Proposition}[section]
\newtheorem{mainproposition}{Main Proposition}[section]
\newtheorem{lemma}{Lemma}[section]
\newtheorem{corollary}{Corollary}[section]
\numberwithin{equation}{section}

\theoremstyle{definition}
\newtheorem{definition}{Definition}[section]

\theoremstyle{remark}
\newtheorem{remark}{Remark}[section]

\def\ba{\begin{array}}
\def\ea{\end{array}}
\def\be{\begin{equation}}
\def\ee{\end{equation}}
\def\bee{\begin{eqnarray}}
\def\beee{\begin{eqnarray*}}
\def\eee{\end{eqnarray}}
\def\eeee{\end{eqnarray*}}
\def\nn{\nonumber}

\title{\bf Dirac-harmonic maps from degenerating spin
surfaces I: the Neveu-Schwarz case \footnotetext{\emph{Date}:
\today.} }
\author{Miaomiao Zhu \thanks{Supported by IMPRS ``Mathematics in
the Sciences'' and the Klaus Tschira Foundation}}
\date{}

\begin{document}
\maketitle

\thispagestyle{empty} \setcounter{page}{1}

\begin{abstract} \vskip 3mm\footnotesize
\noindent We study Dirac-harmonic maps from degenerating spin
surfaces with uniformly bounded energy and show the so-called
generalized energy identity in the case that the domain converges to
a spin surface with only Neveu-Schwarz type nodes. We find condition
that is both necessary and sufficient for the $W^{1,2} \times L^{4}$
modulo bubbles compactness of a sequence of such maps.

\vskip 4.5mm



\end{abstract}

\vskip 0.5cm

\section{Introduction} \label{section 1}

The notion of Dirac-harmonic maps was first introduced in
\cite{CJLW2}. Motivated by the supersymmetric nonlinear sigma model
from quantum field theory \cite{D}, Dirac-harmonic maps are defined
as solutions of a system of harmonic-type equations coupled with
Dirac-type equations. As is done in the theory of minimal surfaces
in Riemannian manifolds and pseudo-holomorphic curves in symplectic
geometry,  construction of geometric invariants from the solution
spaces is expected. This supersymmetric model is introduced in such
a natural way that most fundamental features of two-dimensional
harmonic maps are preserved. Following the approach of Sacks and
Uhlenbeck \cite{SU1}, Chen et al. \cite{CJLW1}, \cite{CJLW2}
developed the ``blow-up" analysis for Dirac-harmonic maps and
established the energy identity for a sequence of Dirac-harmonic
maps with uniformly bounded energy (\cite{CJLW1} for spherical
targets and \cite{Za} for general targets), which gives the $W^{1,2}
\times L^{4}$ modulo bubbles compactness of the solution space for a
fixed spin surface. A natural question then is whether such
compactness is preserved if we allow the domain surface to vary.

To state the problem more precisely, we consider a sequence of
smooth Dirac-harmonic maps
\bee \label{1.1}(\phi_{n},\psi_{n}):
(M_{n},h_{n},c_{n},\mathfrak{S}_{n}) \rightarrow (N,g),\eee
with uniformly bounded energy $E(\phi_{n},\psi_{n},M_{n})\leq
\Lambda <\infty$. Here $(N,g)$ is a compact Riemannian manifold with
metric $g$ and $(M_{n},h_{n},c_{n},\mathfrak{S}_{n})$ is a sequence
of closed hyperbolic Riemann surfaces of genus $g>1$ with hyperbolic
metrics $h_{n}$, compatible complex structures $c_{n}$ and spin
structures $\mathfrak{S}_{n}$.

In this paper, we first prove the energy identity for the sequence
\eqref{1.1} when the domain surface varies in a compact region.
Then, we show the so-called generalized energy identity for the
sequence when the domain surface degenerates to a spin surface with
only Neveu-Schwarz type nodes. The necessary and sufficient
condition for the $W^{1,2} \times L^{4}$ modulo bubbles compactness
of such sequences is a direct consequence of the generalized energy
identity.

Let us consider the simpler case that $(M_{n},h_{n},c_{n})$
converges to a compact hyperbolic Riemann surface $(M,h,c)$ of the
same topological type. Then there exists a sequence of
diffeomorphisms $\tau_{n}: M\rightarrow M_{n}$ such that
$(\tau_{n}^{*}h_{n},\tau_{n}^{*}c_{n})$ converges to $(h,c)$ in
$C^{\infty}$. After passing to a subsequence, we can assume that the
pull-back of $\mathfrak{S}_{n}$ via $\tau_{n}$ is a fixed spin
structure on $M$. Let us denote it by $\mathfrak{S}$. Then, we can
fix the spinor bundle $\Sigma M$ and think of the hyperbolic metrics
$h_{n}$ and the compatible complex structures $c_{n}$ as all living
on the limit surface $M$ and converging in $C^{\infty}$ to $h$ and
$c$, respectively. Let $\nabla_{n}$ be the connection on $\Sigma M$
coming from $h_{n}$ and $\nabla$ the connection on $\Sigma M$ coming
from $h$. Replaced by the pullbacks, we can think of $(\phi_{n},
\psi_{n})$ as a sequence of Dirac-harmonic maps defined on
$(M,h_{n},c_{n},\mathfrak{S})$ with respect to $(c_{n},\nabla_{n})$.
Then we prove the following energy identity for Dirac-harmonic maps
from non-degenerating spin surfaces:

\begin{theorem} \label{thm1.1} Assumptions and notations as above.
Then there exist finitely many blow-up points
$\{x_{1},x_{2},...,x_{I}\}$, finitely many Dirac-harmonic maps
$(\sigma^{i,l},\xi^{i,l}):S^{2}\rightarrow N, i=1,2,...,I;
l=1,2,...,L_{i}$, and a Dirac-harmonic map
$(\phi,\psi):(M,h,\mathfrak{S})\rightarrow N$ such that after
selection of a subsequence, $(\phi_{n},\psi_{n})$ converges to
$(\phi,\psi)$ in $C^{\infty}_{loc}\times C^{\infty}_{loc}$ on $M
\setminus \{x_{1},x_{2},...,x_{I}\}$ and the following hold
\bee
\lim\limits_{n\rightarrow\infty}E(\phi_{n})&=& E(\phi)+\sum\limits_{i=1}^{I}\sum\limits_{l=1}^{L_{i}}E(\sigma^{i,l}),   \\
\lim\limits_{n\rightarrow\infty}E(\psi_{n})&=&
E(\psi)+\sum\limits_{i=1}^{I}\sum\limits_{l=1}^{L_{i}}E(\xi^{i,l}).
\eee
\end{theorem}

To continue the discussions, we recall that the Hopf quadratic
differential associated to a two-dimensional harmonic map plays an
important role in establishing the so-called generalized energy
identity for harmonic maps from degenerating Riemann surfaces
\cite{Z2}. It is observed in \cite{CJLW2} that there is a
generalization of the notion of Hopf differential for a two
dimensional Dirac-harmonic map. Let $(\phi, \psi)$ be a
Dirac-harmonic map defined on a standard cylinder $P=
[t_{1},t_{2}]\times S^{1}$ with flat metric
$ds^{2}=dt^{2}+d\theta^{2}$ and $T(\phi,\psi)(dt+id\theta)^{2}$ the
generalized Hopf differential of $(\phi, \psi)$ on $P$. Then the
following integral
\be \label{1.4}\int_{\{t\}\times S^{1}}T(\phi,\psi) d\theta  \ee
is a complex number which is independent of $t\in[T_{1},T_{2}]$. Let
us denote it by $\alpha=\alpha(\phi, \psi,P)$.

Now we consider the case that $(M_{n},h_{n},c_{n})$ degenerates to a
hyperbolic Riemann surface $(M,h,c)$ by collapsing $p$ ($1\leq p\leq
3g-3$) pairwise disjoint simple closed geodesics $\gamma_{n}^{j}$ of
length $l^{j}_{n}$, $j=1,2,...,p$. For each $j$, the geodesics
$\gamma_{n}^{j}$ degenerate into a pair of punctures
$(\mathcal{E}^{j,1},\mathcal{E}^{j,2})$. Let $P_{n}^{j}$ be the
standard cylindrical collar about $\gamma_{n}^{j}$. Then, we
associate to the sequence $(\phi_{n},\psi_{n},M_{n})$ a sequence of
$p$-tuples $(\alpha_{n}^{1},...,\alpha_{n}^{p})$, where
$\alpha_{n}^{j}:=\alpha(\phi_{n},\psi_{n}, P_{n}^{j}) \in
\mathbb{C}$ are the quantities defined via \eqref{1.4}.

By taking subsequences, we can assume that the pull back of
$\mathfrak{S}_{n}$ via $\tau_{n}$ is a fixed spin structure
$\mathfrak{S}$ on $M$. Note that $M$ has $p$ pairs of punctures. We
require the following additional assumption:
\bee \emph{All punctures of the limit spin surface
$(M,\mathfrak{S})$ are of Neveu-Schwarz type.} \eee
Then $\mathfrak{S}$ extends to some spin structure
$\overline{\mathfrak{S}}$ on $\overline{M}$, where $\overline{M}$ is
the surface obtained by adding a point at each puncture of $M$. As
before, we think of the hyperbolic metrics and the compatible
complex structures $(h_{n},c_{n})$ as all living on the limit
surface $M$ and converging in $C_{loc}^{\infty}$ to $(h,c)$ (c.f.
\cite{Z2}). Thus, $(\phi_{n}, \psi_{n})$ becomes a sequence of
Dirac-harmonic maps defined on $(M,h_{n},c_{n},\mathfrak{S})$. Then
we show the following generalized energy identity for Dirac-harmonic
maps from degenerating spin surfaces:

\begin{theorem} \label{thm1.2} Assumptions and notations as above.
Then there exist finitely many blow-up points
$\{x_{1},x_{2},...,x_{I}\}$ which are away from the punctures
$\{(\mathcal{E}^{j,1},\mathcal{E}^{j,2}), j=1,2,...p\}$ and finitely
many Dirac-harmonic maps

\begin{itemize}
\item[] $(\phi,\psi):(\overline{M},\overline{c},\overline{\mathfrak{S}})\rightarrow N$,
where $(\overline{M},\overline{c},\overline{\mathfrak{S}})$ is the
normalization of $(M,c,\mathfrak{S})$,

\item[]
$(\sigma^{i,l},\xi^{i,l}):S^{2}\rightarrow N, l=1,2,...,L_{i}$, near
the $i$-th blow-up point $x_{i}$,

\item[]
$(\omega^{j,k},\zeta^{j,k}):S^{2}\rightarrow N, k=1,2,...,K_{j}$,
near the $j$-th pair of punctures
$(\mathcal{E}^{j,1},\mathcal{E}^{j,2})$,
\end{itemize}

\noindent such that after selection of a subsequence,
$(\phi_{n},\psi_{n})$ converges to $(\phi,\psi)$ in
$C^{\infty}_{loc}\times C^{\infty}_{loc}$ on $M \setminus
\{x_{1},x_{2},...,x_{I}\}$, and the following holds:
\bee \lim\limits_{n\rightarrow\infty}E(\phi_{n})&=&
E(\phi)+\sum\limits_{i=1}^{I}\sum\limits_{l=1}^{L_{i}}E(\sigma^{i,l})
+\sum_{j=1}^{p}\sum_{k=1}^{K_{j}}E(\omega^{j,k})
+\sum_{j=1}^{p}\lim\limits_{n \rightarrow \infty}|{\rm Re} \alpha
_{n}^{j}|\cdot
\frac{2\pi^{2}}{l_{n}^{j}}, \nn   \\
\\
\lim\limits_{n\rightarrow\infty}E(\psi_{n})&=& E(\psi)
+\sum\limits_{i=1}^{I}\sum\limits_{l=1}^{L_{i}}E(\xi^{i,l})
+\sum_{j=1}^{p}\sum_{k=1}^{K_{j}}E(\zeta^{j,k}). \eee
\end{theorem}

As a corollary, we have

\begin{corollary}  \label{cor1.1} Assumptions and notations as in Theorem \ref{thm1.2}. Then
$(\Phi_{n},\Psi_{n},M_{n})$ subconverge in $W^{1,2}\times L^{4}$
modulo bubbles, i.e., in the limit, the necks contain no energy if
and only if
\bee \liminf\limits_{n \rightarrow \infty}|{\rm Re} \alpha
_{n}^{j}|\cdot \frac{2\pi^{2}}{l_{n}^{j}}=0, \quad j=1,2,...,p. \eee

\end{corollary}

For the asymptotics of the imaginary part of $\alpha_{n}^{j}$, we
have

\begin{proposition} \label{pro1.1} Assumptions and notations as in Theorem \ref{thm1.2}.
 Then
\bee \limsup \limits_{n \rightarrow \infty}|{\rm Im} \alpha
_{n}^{j}|\cdot \frac{2\pi^{2}}{l_{n}^{j}}=0, \quad j=1,2,...,p. \eee
\end{proposition}

We see from the above results that the limits $\liminf\limits_{n
\rightarrow \infty}|{\rm Re} \alpha _{n}^{j}|\cdot
\frac{2\pi^{2}}{l_{n}^{j}}, j=1,2,...,p$ are the obstructions for
$(\phi_{n},\psi_{n},M_{n})$ to subconverge in $W^{1,2} \times L^{4}$
modulo bubbles.

Note that in Theorem \ref{thm1.2}, we made the assumption that all
punctures of the limit surface are of Neveu-Schwarz type. However,
for a general sequence of degenerating spin surfaces, it is possible
that the limit spin surface has Ramond type punctures, along which
the spin structure is trivial. In this case, the corresponding
generalized energy identity is still open.

Now we give a brief outline of the paper. In Sect. 2, we first
recall some preliminary facts about Dirac-harmonic maps from spin
surfaces and then prove Theorem \ref{thm1.1}. In Sect. 3, some
analytic properties of Dirac-harmonic maps from long spin cylinders
are deduced. In Sect. 4, we study Dirac-harmonic maps from
degenerating spin surfaces and show Theorem 1.2.

\vskip 0.2cm

\noindent{\bf Acknowledgements} This paper is part of the author's
Ph.D. thesis \cite{Z1}. He is grateful to his advisor, Prof.
J\"{u}rgen Jost, for guidance and encouragement. He would also like
to thank Prof. Guofang Wang, Prof. Xiaohuan Mo and Guy Buss for
helpful discussions.

\vskip 0.5cm

\section{Notations and preliminaries}
\label{section 2}

In this section, we shall first review some geometric and analytic
aspects of Dirac-harmonic maps and then prove Theorem 1.1.

Let $(M,h,\mathfrak{S})$ be an oriented, compact Riemannian surface
with a fixed spin structure $\mathfrak{S}$ and $P_{\emph{Spin}(2)}
\rightarrow M$ the principal \emph{Spin}(2)-bundle determined by
$\mathfrak{S}$. Let $\Sigma M$ be the spinor bundle over $M$ with a
hermitian metric $\langle\cdot,\cdot\rangle_{\Sigma M}$. The
Levi-Civita connection $\nabla^{TM}$ on $TM$ with respect to $h$
gives rise to a connection-1-form
$\{\omega_{\alpha\beta}\}_{\alpha,\beta=1}^{2}$ on
$P_{\emph{Spin}(2)}$, and this in turn defines a spin connection
$\nabla^{\Sigma M}$ on $\Sigma M$ that is compatible with
$\langle\cdot,\cdot\rangle_{\Sigma M}$. For simplicity of notation,
we denote $\nabla^{\Sigma M}$ by $\nabla$. The Dirac operator
$\slashed{\partial}$ is locally given by $ \slashed{\partial} \psi:=
e_{1} \cdot \nabla _{e_{1}} \psi+ e_{2} \cdot \nabla _{e_{2}}\psi$
for a local orthonormal frame $\{e_{1},e_{2}\}$ of $TM$ and
$\psi\in\Sigma M$. We refer to \cite{LM}, \cite{Gi}, and \cite{F}
for more background material on spin structures and Dirac operators
and to \cite{J} for general Riemannian geometrical notations.

Let $(N,g)$ be a compact Riemannian manifold of dimension $d\geq 2$
and $\phi$ a smooth map from $M$ to $N$. By $\phi^{-1}TN$, we denote
the pull-back bundle of $TN$ via $\phi$. Consider the twisted bundle
$\Sigma M\otimes \phi^{-1}TN$ with a metric
$\langle\cdot,\cdot\rangle_{\Sigma M\otimes \phi^{-1}TN}$ induced
from the metrics on $\Sigma M$ and $\phi^{-1}TN$. There is a natural
connection $\widetilde{\nabla}$ on $\Sigma M\otimes \phi^{-1}TN$
induced from those on $\Sigma M$ and $\phi^{-1}TN$, namely,
\bee \widetilde{\nabla}:=\nabla^{\Sigma M}\otimes 1 + 1 \otimes
\nabla^{\phi^{-1}TN}. \nn \eee
The section $\psi\in \Gamma(\Sigma M \otimes \phi^{-1}TN)$ is
written in local coordinates as
$\psi=\psi^{i}\otimes\partial_{y^{i}}(\phi)$, where $\psi^{i}\in
\Sigma M$ and $\{\partial_{y^{i}} \}$ is a local basis on $N$. The
\emph{Dirac operator along the map} $\phi$ is defined by
\bee \label{2.7}\slashed{D}\psi := e_{\alpha}\cdot
\widetilde{\nabla}_{e_{\alpha}}\psi,  \nn \eee
where $\psi \in \Gamma(\Sigma M \otimes \phi^{-1}TN)$. Here and in
the sequel, we apply the Einstein summation convention. Set
\bee \chi(M,N):=\left\{(\phi,\psi)|\phi \in C^{\infty}(M,N), \psi
\in \Gamma(\Sigma M \otimes \phi^{-1}TN)\right\} \nn \eee
and consider the following functional defined on $\chi(M,N)$:
\bee L(\phi,\psi) &:=&
\int\limits_{M}\left(|d\phi|^{2}+\langle\psi,\slashed{D}
\psi\rangle_{\Sigma M\otimes \phi^{-1}TN}\right)d vol(h)          \nn  \\
&=& \int\limits_{M}\left(g_{ij}(\phi)h^{\alpha\beta}\frac{\partial
\phi^{i}}{\partial x_{\alpha}}\frac{\partial \phi^{j}}{\partial
x_{\beta}}+g_{ij}(\phi)\langle\psi^{i},\slashed{D}
\psi^{j}\rangle_{\Sigma M}\right)\sqrt{{\rm
det}(h_{\alpha\beta})}dx^{1}dx^{2}. \nn
 \eee
By a straightforward computation (see \cite{CJLW2}), we get the
Euler-Lagrange equations of $L$:
\bee \label{2.1} \tau(\phi) &=& \mathcal{R}(\phi,\psi),     \\
\label{2.2} \slashed{D}\psi &=& 0, \eee
where $\tau(\phi)\in\Gamma(\phi^{-1}TN)$ is the tension field of
$\phi$ and $\mathcal{R}(\phi,\psi)\in\Gamma(\phi^{-1}TN)$ is defined
by
\bee \mathcal{R}(\phi,\psi)(x):=\frac{1}{2}
R^{m}_{lij}(\phi(x))\langle\psi^{i},\nabla\phi^{l}\cdot\psi^{j}\rangle\partial_{
y^{m}}(\phi(x)). \nn \eee
Here $R^{m}_{lij}$ are the components of the curvature tensor of $g$
and $\nabla \phi^{l} \cdot \psi^{j}$ denotes the Clifford
mutiplication of the vector field $\nabla
\phi^{l}:=\phi^{l}_{\alpha}e_{\alpha}$ with the spinor $\psi^{j}$.

Solutions $(\phi,\psi)$ of \eqref{2.1}, \eqref{2.2} are called
\textit{Dirac-harmonic maps} from $M$ to $N$. Thus, a Dirac-harmonic
map is a map coupled with a spinor field with values in the
pull-back tangent bundle. For nontrivial examples, see \cite{CJLW2}.

By the Nash-Moser embedding theorem, we embed $N$ into some
$\mathbb{R}^{K}$. Let $A(\cdot ,\cdot)$ be the second fundamental
form of $N$ in $\mathbb{R}^{K}$ and $P(\cdot ;\cdot)$ the shape
operator, satisfying $\langle P(\xi;X),Y\rangle=\langle
A(X,Y),\xi\rangle$ for any $X,Y\in \Gamma(TN),
\xi\in\Gamma(T^{\perp}N)$, where $T^{\perp}N$ is the normal bundle.
Set
\bee
\mathcal{A}(d\phi(e_{\alpha}),e_{\alpha}\cdot\psi)&:=&\phi^{i}_{\alpha}e_{\alpha}\cdot\psi^{j}\otimes
A(\partial_{y^{i}},\partial_{y^{j}}), \nn  \\
\mathcal{P}(\mathcal{A}(d\phi(e_{\alpha}),e_{\alpha}\cdot\psi);\psi)
&:=&P(A(\partial_{y^{l}},\partial_{y^{j}});\partial_{y^{i}})\langle\psi^{i},e_{\alpha}\cdot\psi^{j}\rangle_{\Sigma
M}\phi_{\alpha}^{l}.  \nn  \eee
Then equations \eqref{2.1} and \eqref{2.2} become
\bee \label{2.3} -\Delta\phi &=& A(d\phi,d\phi)+{\rm
Re}\ \mathcal{P}(\mathcal{A}(d\phi(e_{\alpha}),e_{\alpha}\cdot\psi);\psi)  \\
\label{2.4}   \slashed{\partial}\psi &=&
\mathcal{A}(d\phi(e_{\alpha}),e_{\alpha}\cdot\psi).\eee
Here, $\phi$ is a map from $M$ to $\mathbb{R}^{K}$ with
\bee \label{2.5} \phi(x)\in N \eee
for any $x\in M$, and the spinor field $\psi$ along the map $\phi$
is a $K$-tuple of spinors $(\psi^{1},\psi^{2},...,\psi^{K})$
satisfying
\bee \label{2.6} \sum_{i}\nu_{i}\psi^{i}=0, \quad {\rm for\ any\
normal\ vector}\ \nu=\sum_{i=1}^{K}\nu_{i}E_{i}\ {\rm at}\ \phi(x),
\eee
where $\{E_{i},i=1,2,...,K\}$ is the standard basis of
$\mathbb{R}^{K}$.

Set
\bee \chi_{1,4/3}^{1,2}(M,N):=\{(\phi,\psi)\in W^{1,2}\times
W^{1,4/3}\ {\rm with}\ \eqref{2.5}\ {\rm and}\ \eqref{2.6}\ {\rm
a.e.} \}. \nn  \eee
Then the functional $L(\phi,\psi)$ is well-defined for
$(\phi,\psi)\in \chi_{1,4/3}^{1,2}(M,N)$. A critical point
$(\phi,\psi)$ of the functional $L$ in $\chi_{1,4/3}^{1,2}(M,N)$ is
called a \emph{weakly Dirac-harmonic map} from $M$ to $N$. When the
target $N$ is the standard sphere $\mathbb{S}^{d}$, a weakly
Dirac-harmonic map is smooth \cite{CJLW1}.

\vskip 0.2cm

Let $\Omega$ be a domain of $M$. The energy of $(\phi,\psi)$ on
$\Omega$ is defined by
\bee E(\phi,\psi,\Omega):=\int_{\Omega}(|d\phi|^{2}+|\psi|^{4}).
  \nn \eee
The energy of $\phi$ on $\Omega$ is
\bee E(\phi,\Omega):=\int_{\Omega}|d\phi|^{2} \nn  \eee
and the energy of $\psi$ on $\Omega$ is
\bee E(\psi,\Omega):=\int_{\Omega}|\psi|^{4}.   \nn  \eee

For a two-dimensional harmonic map, there are two important
geometric properties, namely, the conformal invariance and the
existence of the Hopf quadratic differential. The following two
propositions, proved in \cite{CJLW2}, show that these two properties
are preserved in the case of Dirac-harmonic maps.

\begin{proposition} \label{pro2.1} The functional $L(\phi,\psi)$ and the energy
$E(\phi,\psi)$ are conformally invariant. Namely, for any conformal
diffeomorphism $ f: M \rightarrow M$, set
\bee \widetilde{\phi}=\phi\circ f, \quad
\widetilde{\psi}=e^{-\frac{\sigma}{2}}\psi\circ f.  \nn \eee
Then $L(\phi,\psi)=L(\widetilde{\phi},\widetilde{\psi}),
E(\phi,\psi)=E(\widetilde{\phi},\widetilde{\psi})$. Here
$e^{-\sigma}$ is the conformal factor of the conformal map $f$.
\end{proposition}

\begin{remark} In fact, the following terms are all conformally
invariant:
\bee \int |d \phi|^{2}d vol(h), \int
\langle\psi,\slashed{D}\psi\rangle d vol(h), \int |\psi|^{4}d
vol(h). \nn  \eee
\end{remark}

Let $(\phi,\psi)$ be a Dirac-harmonic map from $(M,h)$. Let $\Omega
\subset M$ be a small domain, and take a local isothermal coordinate
$z=x+iy$ on $\Omega$ such that $h=\rho|dz|^{2}$. Define
\be \label{2.35} T(\phi,\psi)(z)dz^{2}=
\left\{(|\phi_{x}|^2-|\phi_{y}|^{2}-2i\phi_{x}\cdot \phi_{y})+({\rm
Re}\langle\psi,\partial_{x}\cdot\psi_{x} \rangle-i{\rm Re}
\langle\psi,\partial_{x}\cdot\psi_{y}\rangle)\right\}dz^{2}.  \nn
\ee
Here $\partial_{x}=\frac{\partial}{\partial_{x}},
\partial_{y}=\frac{\partial}{\partial_{y}},  \psi_{x}=
\widetilde{\nabla}_{\partial_{x}}\psi,  \psi_{y}=
\widetilde{\nabla}_{\partial_{y}}\psi$. Then we have

\begin{proposition} \label{pro2.2} $T(\phi,\psi)(z)dz^{2}$
is a holomorphic quadratic differential.
\end{proposition}

Now we turn to some analytic aspects of Dirac-harmonic maps. In
\cite{CJLW1} and \cite{CJLW2}, several basic properties of
Dirac-harmonic maps which play an important role in the ``bubbling"
process were established. They can be considered as a generalization
of the corresponding properties of harmonic maps. For the sake of
completeness, we present them here.

\begin{proposition} \label{pro2.3} Let $(M,h)$ be a Riemann surface with a fixed spin structure and
$(N,g)$ be a compact Riemannian manifold of dimension $d$. Then
there is a small constant $\epsilon_{0} >0$ such that if
$(\phi,\psi): M \rightarrow N$ is a smooth Dirac-harmonic map
satisfying
\bee \label{2.36} \int_{M}(|d\phi|^{2}+|\psi|^{4})<\epsilon_{0},
 \nn \eee
then $\phi$ is constant and consequently $\psi$ is a $d$-tuple of
harmonic spinors.
\end{proposition}

\begin{theorem} \label{thm2.1}($\epsilon$-regularity theorem) There is a small constant $\epsilon_{0}>0$
such that if $(\phi,\psi):(D,\delta_{\alpha\beta})\rightarrow
(N,g_{ij})$ is a smooth Dirac-harmonic map satisfying
\bee \int_{D}(|d\phi|^{2}+|\psi|^{4})<\epsilon_{0}, \nn \eee
then
\bee \|d\phi\|_{\widetilde{D},1,p} &\leq&
C(\widetilde{D},p)\|d\phi\|_{D,0,2},          \nn \\
\|\nabla\psi\|_{\widetilde{D},1,p} &\leq&
C(\widetilde{D},p)\|\psi\|_{D,0,4},              \nn \\
\|\nabla\psi\|_{L^{\infty}(\widetilde{D})} &\leq&
C(\widetilde{D})\|\psi\|_{D,0,4},                 \nn  \\
\|\psi\|_{L^{\infty}(\widetilde{D})} &\leq&
C(\widetilde{D})\|\psi\|_{D,0,4}, \nn  \eee
$\forall\widetilde{D}\subset \subset D,p>1$, where
$C(\widetilde{D},p)>1$ is a constant depending only on
$\widetilde{D},p$, and the geometry of $N$.
\end{theorem}

Before we state the theorem on the removability of isolated
singularities, let us recall some facts about the spin structures on
surfaces (c.f. \cite{AB}, Sect. 2).

Let $(M,h)$ be an oriented Riemannian surface and $P_{SO(2)}$ its
oriented orthonormal frame bundle. Let $\gamma : S^{1}\rightarrow M$
be an immersion. Then the unit tangent vector field of $\gamma$
together with the corresponding unit normal vector field forms a
section of $P_{SO(2)}$ along $\gamma$. A spin structure of $M$ is
said to be \emph{trivial} along $\gamma$ if this section lifts to a
closed curve in $P_{\emph{Spin}(2)}$; otherwise, it is said to be
\emph{nontrivial} along $\gamma$. It should be remarked that this
notion is invariant under deformations of $\gamma$ within the same
homotopy class of immersions and hence can be used to specify the
two different spin structures on an annulus or a cylinder. There are
various equivalent definitions of the triviality of a spin structure
along a cylindrical end, see for instance \cite{B}.

Now we consider a punctured disk $D\setminus \{0\}$ with the spin
structure being nontrivial along $\partial D$. Note that this spin
structure extends to the unique spin structure on $D$ (c.f.
\cite{AB}, Sect. 2). Then we have

\begin{theorem} \label{thm2.2} (Removable singularity theorem) Let
$(\phi,\psi)$ be a solution of \eqref{2.1} and \eqref{2.2} which is
$C^{\infty}$ on $D\setminus \{0\}$. If $(\phi,\psi)$ has finite
energy, then $(\phi,\psi)$ extends to a $C^{\infty}$ solution on
$D$.
\end{theorem}

\begin{remark} Here, the singularity
$\{0\}$ is said to be of \emph{Neveu-Schwarz} type (see \cite{JKV}
for an algebraic geometric description). However, there is another
spin structure on $D\setminus \{0\}$ that cannot be extended to the
unique spin structure on $D$, and the singularity $\{0\}$ then is
said to be of \emph{Ramond} type \cite{JKV}. We do not know whether
an analogous theorem for the Ramond type singularities also holds.
\end{remark}

For the proofs of Proposition \ref{pro2.3}, Theorem \ref{thm2.1} and
Theorem \ref{thm2.2}, see \cite{CJLW2}.

Applying the geometric and analytic properties of Dirac-harmonic
maps developed before, Chen et al. \cite{CJLW1} and Zhao \cite{Za}
studied the compactness of a sequence of smooth Dirac-harmonic maps
from a fixed domain and proved the following energy identity
theorem.

\begin{theorem} \label{thm2.3} Let $(\phi_{n},\psi_{n}):(M,h,\mathfrak{S})\rightarrow N$
 be a sequence of smooth Dirac-harmonic maps with uniformly bounded energy
$E(\phi_{n},\psi_{n})\leq \Lambda < +\infty.$ Then there exist
finitely many blow-up points $\{x_{1},x_{2},...,x_{I}\}$, finitely
many Dirac-harmonic maps $(\sigma^{i,l},\xi^{i,l}):S^{2}\rightarrow
N, i=1,2,...,I; l=1,2,...,L_{i}$, and a Dirac-harmonic map
$(\phi,\psi):(M,h,\mathfrak{S})\rightarrow N$ such that, after
selection of a subsequence,
$(\phi_{n},\psi_{n})\rightarrow(\phi,\psi)$ in
$C^{\infty}_{loc}\times C^{\infty}_{loc}$ on
$M\setminus\{x_{1},x_{2},...,x_{I}\}$, and the following holds:
\bee
\lim\limits_{n\rightarrow\infty}E(\phi_{n})&=& E(\phi)+\sum\limits_{i=1}^{I}\sum\limits_{l=1}^{L_{i}}E(\sigma^{i,l}),   \\
\lim\limits_{n\rightarrow\infty}E(\psi_{n})&=&
E(\psi)+\sum\limits_{i=1}^{I}\sum\limits_{l=1}^{L_{i}}E(\xi^{i,l}).
\eee
\end{theorem}
\begin{remark} \label{rem2.3}
When the domain is fixed, the ``bubbling" procedure corresponds to
collapsing homotopically trivial simple closed curves on the domain
surface. During this process, some necks joining one bubble to the
next appear, as in the case of harmonic maps. On the one hand, by
applying the standard ``blow-up" analysis on cylinders (c.f. Theorem
3.6 in \cite{CJLW1}), we can obtain Dirac-harmonic maps from
$\mathbb{R}\times S^{1}$. On the other hand, since any spin
structure on a surface is nontrivial along any homotopically trivial
simple closed curve (c.f. \cite{AB}, Sect. 2), the induced spin
structure on each $\mathbb{R}\times S^{1}$ is nontrivial. Also, the
induced spin structures on the domain cylinders of the necks are
nontrivial. Note that the nontrivial spin structure on
$\mathbb{R}\times S^{1}$ can be conformally compactified to the
unique spin structure on $S^{2}$. Thus, one can apply the conformal
invariance of Dirac-harmonic maps and the removable singularity
theorem to obtain Dirac-harmonic maps from $S^{2}$. The
nontriviality of the spin structures along the domain cylinders is
crucial here.
\end{remark}

It is interesting to ask what happens when the domain of the
Dirac-harmonic maps $(\phi_{n},\psi_{n})$ varies. By Riemann surface
theory, we can fix the topological type of the surface and let the
complex structure of the surface vary with $n$. The conformal
invariance of Dirac-harmonic maps allows us to take a particular
metric within the same conformal class. To do this, we consider
closed Riemann surfaces of genus $g>1$. It follows from the
uniformization theorem that any such surface acquires a complete
hyperbolic metric that is unique in the conformal class determined
by the complex structure. Thus, we have the following data
associated to a spin surface:
\bee (M,h,c,\mathfrak{S}). \nn \eee
Here, $c$ is a complex structure, $h$ is the hyperbolic metric
compatible with $c$ and $\mathfrak{S}$ is a spin structure.

Now, we consider a sequence of smooth Dirac-harmonic maps
\bee (\phi_{n},\psi_{n}): (M_{n},h_{n},c_{n},\mathfrak{S}_{n})
\rightarrow N    \nn \eee
with uniformly bounded energy $E(\phi_{n},\psi_{n},M_{n})\leq
\Lambda <\infty.$

Let us consider the simpler case that $(M_{n},h_{n},c_{n})$
converges to a closed hyperbolic Riemann surface $(M,h,c)$ of the
same topological type. Then there exists a sequence of
diffeomorphisms $\tau_{n}: M\rightarrow M_{n}$ such that
$(\tau_{n}^{*}h_{n},\tau_{n}^{*}c_{n})$ converges to $(h,c)$ in
$C^{\infty}$ (c.f. \cite{Z2}). We need to consider the change of the
spin structure involved. In general, a diffeomorphism between two
spin surfaces may not preserve the spin structures. However, for a
closed Riemann surface of genus $g$, there are exactly $2^{2g}$
topologically different equivalence classes of spin structures
\cite{LM}. Hence, after passing to a subsequence, we can assume that
the diffeomorphisms $\tau_{n}$ are compatible with the spin
structures $\mathfrak{S}_{n}$, namely, the pull back of
$\mathfrak{S}_{n}$ via $\tau_{n}$ is a fixed spin structure on $M$.
We denote it by $\mathfrak{S}$. Recall that to fix a spin structure
means to fix the equivalence class of a spin structure, thus the
corresponding principle $\emph{Spin}(2)$-bundles can be naturally
identified with each other via bundle isomorphisms (c.f. \cite{BG}
or \cite{M}). Likewise, the corresponding associated spinor bundles
can also be identified with each other. As explained in \cite{L}, we
think of the principle $\emph{Spin}(2)$-bundle $P_{\emph{Spin}(2)}$
as a topological fiber bundle and $\Sigma M$ as the associated
bundle with a natural hermitian metric $\langle \cdot,\cdot
\rangle_{\Sigma M}$. They are independent of the metric $h$ chosen,
as long as the spin structure is fixed. The metric $h$ enters in
defining the connection-1-form $\{\omega_{\alpha\beta}\}$ and hence
the spin connection $\nabla^{\Sigma M}$. Thus, we can fix the spinor
bundle $\Sigma M$ and think of the hyperbolic metrics $h_{n}$ and
the compatible complex structures $c_{n}$ as all living on the limit
surface $M$ and converging in $C^{\infty}$ to $h$ and $c$,
respectively. Let $\nabla_{n}$ be the connection on $\Sigma M$
coming from $h_{n}$ and $\nabla$ be the connection on $\Sigma M$
coming from $h$. Replaced by the pullbacks, we think of $(\phi_{n},
\psi_{n})\in C^{\infty}(M,N)\times C^{\infty}(\Sigma M\otimes
\mathbb{R}^{K})$ as a sequence of Dirac-harmonic maps defined on
$(M,h_{n},c_{n},\mathfrak{S})$ with respect to $(c_{n},\nabla_{n})$.

\vskip 0.2cm

\noindent \emph{Proof of Theorem \ref{thm1.1}}. As $n\rightarrow
\infty$, $(h_{n},c_{n})$ converges in $C^{\infty}$ to $(h,c)$.
Hence, all geometric data associated to $(h_{n},c_{n})$ converge in
$C^{\infty}$ to those associated to $(h,c)$. In particular, the
tensor $\nabla_{n}-\nabla \in {\rm End}(\Sigma M,\Sigma M\otimes
T^{*}M)$ converges to zero in $C^{\infty}$ and the energy functional
corresponding to $(c_{n},\nabla_{n})$ is uniformly equivalent to the
one corresponding to $(c,\nabla)$. By the uniform energy bound
$E(\phi_{n},\psi_{n},M_{n})\leq \Lambda$, we can assume that
$(\phi_{n},\psi_{n})$ weakly converges to some $(\phi,\psi)$ in
$W^{1,2}(M,N)\times L^{4}(\Sigma M\otimes \mathbb{R}^{K})$ with
respect to $(c,\nabla)$. Note that all estimates in Proposition
\ref{pro2.3}, Theorem \ref{thm2.1} and Theorem \ref{thm2.3} are
uniform for the metrics $h_{n}$ and the complex structures $c_{n}$.
Hence, by the standard covering argument (c.f. \cite{SU1} and
Theorem 2.3 in \cite{SU2}), there exist finitely many points
$\{x_{1},x_{2},...,x_{I}\}$ in $M$ such that $(\phi_{n},\psi_{n})$
subconverges in $C^{\infty}_{loc}\times C^{\infty}_{loc}$ to
$(\phi,\psi)$ on $M \setminus \{x_{1},x_{2},...,x_{I}\}$. By the
smoothness of $(\phi_{n},\psi_{n})$, we know that $(\phi,\psi)\in
W^{1,2}(M,N)\times L^{4}(\Sigma M\otimes \mathbb{R}^{K})$ is
actually a smooth Dirac-harmonic map defined on $M \setminus
\{x_{1},x_{2},...,x_{I}\}$ with respect to $(c,\nabla)$. By the
removable singularity theorem, $(\phi,\psi)$ extends to a smooth
Dirac-harmonic map from $M$ with respect to $(c,\nabla)$. The rest
of the proof of the theorem is almost immediate from applying the
``blow-up" process to capture the energy concentration at the
isolated singularities, which is analogous to the proof of Theorem
\ref{thm2.3} (see \cite{CJLW1}, \cite{Za}), since all estimates are
uniform for $(h_{n},c_{n})$. \hfill $\square$

\vskip 0.5cm

\section{Dirac-harmonic maps from spin cylinders}
\label{section 3}

In this section, we establish a series of analytic properties of
Dirac-harmonic maps from spin cylinders.

For later use, we introduce a conformal transformation between an
annulus and a cylinder. Let $(r,\theta)$ be the polar coordinates of
$\mathbb{R}^{2}$ centered at 0 and
$h_{eucl}=dr^{2}+r^{2}d\theta^{2}$ the Euclidean metric on
$\mathbb{R}^{2}$. Consider a map $f: \mathbb{R}^{1}\times
S^{1}\rightarrow \mathbb{R}^{2}$ given by
\be \label{3.1} r=e^{-t},\quad \theta=\theta,\quad (t,\theta)\in
\mathbb{R}^{1}\times S^{1}.\ee
Let us equip $\mathbb{R}^{1}\times S^{1}$ with the metric $
ds^{2}=dt^{2}+d\theta^{2}.$ Then it is easy to verify that
\be f^{*}h_{eucl}=e^{-2t}ds^{2}. \nn \ee
Thus $f:\mathbb{R}^{1}\times S^{1}\rightarrow \mathbb{R}^{2}$ is a
conformal transformation. Given $r_{1}>r_{2}$, then, the annulus
$A_{r_{1},r_{2}}:=\{re^{i\theta}| r_{2}\leq r\leq r_{1}\}$ is mapped
to the cylinder $P_{t_{1},t_{2}}:=[t_{1},t_{2}]\times S^{1}$, where
$t_{i}=-\log r_{i},i=1,2.$

Let $(\phi,\psi)$ be a Dirac-harmonic map defined on the annulus
$A_{r_{1},r_{2}}\subset\mathbb{R}^{2}$. Set
\be \Phi:=f^{*}\phi, \quad  \Psi:=e^{-\frac{t}{2}}f^{*}\psi. \nn \ee
Then by the conformal invariance of Dirac-harmonic maps,
$(\Phi,\Psi)$ is a Dirac-harmonic map defined on the cylinder
$P_{t_{1},t_{2}}\subset\mathbb{R}^{1}\times S^{1}$.

By $P_{T_{1},T_{2}}=[T_{1},T_{2}]\times S^{1}$, we denote a cylinder
with metric $ds^{2}=dt^{2}+d\theta^{2}$ and with the spin structure
being nontrivial along the boundary curves.

The following lemma is a cylindrical version of Lemma 3.2 in
\cite{Za}:

\begin{lemma} \label{lem3.1} Let $(\Phi,\Psi)\in
C^{\infty}(P_{T_{1},T_{2}},N)$ be a Dirac-harmonic map, where
$T_{2}-1>T_{1}>0$. Then we have
\bee (\int\limits_{P_{T_{1},T_{2}-1}}|\Psi|^{4})^{\frac{1}{4}}
&\leq& C_{0}(\int\limits_{P_{T_{1},T_{2}}}|d\Phi|^{2})^{\frac{1}{2}}
(\int\limits_{P_{T_{1},T_{2}}}|\Psi|^{4})^{\frac{1}{4}}+ C(\int\limits_{P_{T_{2}-1,T_{2}}}|\Psi|^{4})^{\frac{1}{4}}      \nn  \\
&&{}+
C(\int\limits_{S_{1}}|\nabla\Psi|^{\frac{4}{3}})^{\frac{3}{4}}+C(\int\limits_{S_{1}}
|\Psi|^{4})^{\frac{1}{4}}, \eee
\bee
(\int\limits_{P_{T_{1},T_{2}-1}}|\nabla\Psi|^{\frac{4}{3}})^{\frac{3}{4}}
&\leq& C_{0}(\int\limits_{P_{T_{1},T_{2}}}|d\Phi|^{2})^{\frac{1}{2}}(\int\limits_{P_{T_{1},T_{2}}}|\Psi|^{4})^{\frac{1}{4}}
+ C(\int\limits_{P_{T_{2}-1,T_{2}}}|\Psi|^{4})^{\frac{1}{4}}                                              \nn  \\
&&{}+
C(\int\limits_{S_{1}}|\nabla\Psi|^{\frac{4}{3}})^{\frac{3}{4}}+C(\int\limits_{S_{1}}|\Psi|^{4})^{\frac{1}{4}},\eee
where $S_{i}=\{T_{i}\}\times S^{1}, i=1,2$, and $C_{0}, C$ are
constants that do not depend on $T_{1}$ and $T_{2}$.
\end{lemma}

\noindent \emph{Proof}. The result follows from applying the
conformal transformation \eqref{3.1} to Lemma 3.2 in \cite{Za}.
\hfill $\square$

Moreover, let $q(t)$ be an $\mathbb{R}^{K}$-valued linear function
on $[T_{1},T_{2}]$ such that $q(T_{i})$ equals the mean value of
$\Phi$ over $S_{i}, i=1,2$. By employing the technique used by Sacks
and Uhlenbeck in \cite{SU1}, we have the following lemma, which is
part of Lemma 3.3 in \cite{Za}.

\begin{lemma} \label{lem3.2} Let $(\Phi,\Psi)$ be a
Dirac-harmonic map defined on $P_{T_{1},T_{2}}$, where
$T_{2}>T_{1}>0$. Then we have
\bee \int\limits_{P_{T_{1},T_{2}}}|\Phi_{\theta}|^{2} &\leq& C
\underset{P_{T_{1},T_{2}}}{\rm
sup}|\Phi-q|\int\limits_{P_{T_{1},T_{2}}}|d\Phi|^{2}
+C\underset{P_{T_{1},T_{2}}}{\rm sup}|\Phi-q|\int\limits_{P_{T_{1},T_{2}}}|\Psi|^{4}      \nn  \\
&&{}+\int\limits_{S_{1}}-\int\limits_{S_{2}}(\Phi-q)\Phi_{t}d\theta.
\eee
Here $C$ is a constant that only depends on $N$, not on $T_{1}$ and
$T_{2}$.
\end{lemma}

For a proof, see \cite{Za}. Note that here we only estimate the
vertical energy.

Inspired by the proof of Theorem 3.5 in \cite{Za}, we give the
following lemma.
\begin{lemma}\label{lem3.3} There exists $\epsilon_{1}>0$ such
that if $(\Phi,\Psi)$ is a Dirac-harmonic map defined on
$P_{T_{1}-1,T_{2}+1}$ and
\be \label{3.5}
\int\limits_{P_{T_{1}-1,T_{2}+1}}|d\Phi|^{2}+|\Psi|^{4} \leq \Lambda
< \infty, \ee
\be \label{3.6}\omega:=\underset{t\in[T_{1}-1,T_{2}]}{\rm
sup}\int\limits_{[t,t+1]\times S^{1}}|d\Phi|^{2}+|\Psi|^{4} \leq
\epsilon_{1}, \ee
then
\be \label{3.7} \int\limits_{P_{T_{1},T_{2}}}|\Phi_{\theta}|^{2}
+\int\limits_{P_{T_{1},T_{2}}}|\Psi|^{4}+\int\limits_{P_{T_{1},T_{2}}}|\nabla\Psi|^{\frac{4}{3}}
\leq C(\Lambda)\omega^{\frac{1}{3}}. \ee
Here, $C(\Lambda)$ is a constant depending only on $\Lambda$, but
not on $T_{1},T_{2}$.
\end{lemma}

\noindent\emph{Proof}. Let $\epsilon_{1}={\rm
min}\{\epsilon_{0},\frac{1}{8C_{0}^{2}},1\}$, where $\epsilon_{0}$
is the constant in the $\epsilon$-regularity theorem and $C_{0}$ is
the constant in Lemma \ref{lem3.1}. Then the assumption \eqref{3.6}
implies
\bee \label{3.8} \underset{t\in[T_{1}-1,T_{2}]}{\rm
sup}\int\limits_{[t,t+1]\times S^{1}}|d\Phi|^{2}+|\Psi|^{4} \leq
\epsilon_{1}\leq \frac{1}{8C_{0}^{2}}. \eee
Since $ \mu(t):=\int_{[T_{1},t]\times S^{1}}|d\Phi|^{2}$ is a
continuous and nondecreasing function on $[T_{1}, T_{2}]$ and the
energy of $\Phi$ over $P_{T_{1}-1,T_{2}+1}$ is bounded by $\Lambda$,
we can separate $P_{T_{1},T_{2}}$ into finitely many parts (c.f.
\cite{Za}, p. 134 or similar arguments in \cite{Ye}, p. 689)
\be P_{T_{1},T_{2}}=\bigcup\limits_{n=1}^{N_{0}}P^{n},
P^{n}:=[T^{n-1},T^{n}]\times S^{1}, T^{0}=T_{1}, T^{N_{0}}=T_{2} \nn
\ee
such that $N_{0}$ is an integer no larger then
$[8C_{0}^{2}\Lambda]+1$, and
\be \label{3.9} E(\Phi;P^{n})\leq \frac{1}{4C_{0}^{2}}, \quad
n=1,2,...,N_{0}. \ee
On each part $P^{n}$, by Lemma \ref{lem3.1}, we have
\bee \label{3.10} (\int\limits_{P^{n}}|\Psi|^{4})^{\frac{1}{4}}
&\leq& C_{0}(\int\limits_{[T^{n-1},T^{n}+1]\times
S^{1}}|d\Phi|^{2})^{\frac{1}{2}}
(\int\limits_{[T^{n-1},T^{n}+1]\times S^{1}}|\Psi|^{4})^{\frac{1}{4}}  \nn  \\
&&{}+ C(\int\limits_{[T^{n},T^{n}+1]\times S^{1}}|\Psi|^{4})^{\frac{1}{4}}         \nn  \\
&&{}+ C(\int\limits_{T^{n-1}\times
S^{1}}|\nabla\Psi|^{\frac{4}{3}})^{\frac{3}{4}}+C(\int\limits_{T^{n-1}\times
S^{1}}|\Psi|^{4})^{\frac{1}{4}}.    \eee
Note that $[T^{n-1},T^{n}+1]\times S^{1}= P^{n} \cup
([T^{n},T^{n}+1]\times S^{1})$. Hence, by the following
inequalities:
\bee  (a+b)^{\frac{1}{2}}\leq
(a^{\frac{1}{2}}+b^{\frac{1}{2}}),\quad (a+b)^{\frac{1}{4}}\leq
(a^{\frac{1}{4}}+b^{\frac{1}{4}}), \quad \forall a,b \geq 0. \nn
\eee
and H\"{o}lder's inequality, it is easy to verify that
\bee \label{3.11}(\int\limits_{P^{n}}|\Psi|^{4})^{\frac{1}{4}}
&\leq&
C_{0}(\int\limits_{P^{n}}|d\Phi|^{2})^{\frac{1}{2}}(\int\limits_{P^{n}}|\Psi|^{4})^{\frac{1}{4}}  \nn  \\
&&{}+ C_{0}(\int\limits_{P^{n}}|d\Phi|^{2})^{\frac{1}{2}}(\int\limits_{[T^{n},T^{n}+1]\times S^{1}}|\Psi|^{4})^{\frac{1}{4}}  \nn  \\
&&{}+ C_{0}(\int\limits_{[T^{n},T^{n}+1]\times S^{1}}|d\Phi|^{2})^{\frac{1}{2}}(\int\limits_{P^{n}}|\Psi|^{4})^{\frac{1}{4}}  \nn  \\
&&{}+ C_{0}(\int\limits_{[T^{n},T^{n}+1]\times S^{1}}|d\Phi|^{2})^{\frac{1}{2}}
(\int\limits_{[T^{n},T^{n}+1]\times S^{1}}|\Psi|^{4})^{\frac{1}{4}}  \nn  \\
&&{}+ C(\int\limits_{[T^{n},T^{n}+1]\times S^{1}}|\Psi|^{4})^{\frac{1}{4}}                                          \nn  \\
&&{}+C(\int\limits_{T^{n-1}\times
S^{1}}|\nabla\Psi|^{\frac{4}{3}})^{\frac{3}{4}}+C(\int\limits_{T^{n-1}\times
S^{1}} |\Psi|^{4})^{\frac{1}{4}}. \eee
From \eqref{3.5}, \eqref{3.6}, and \eqref{3.9}, we can rewrite
\eqref{3.11} as follows:
\bee \label{3.12} (\int\limits_{P^{n}}|\Psi|^{4})^{\frac{1}{4}}
&\leq& C(\int\limits_{[T^{n},T^{n}+1]\times
S^{1}}|\Psi|^{4})^{\frac{1}{4}}
+ C(\int\limits_{[T^{n},T^{n}+1]\times S^{1}}|d\Phi|^{2})^{\frac{1}{2}}  \nn  \\
&&{}+ C(\int\limits_{T^{n-1}\times
S^{1}}|\nabla\Psi|^{\frac{4}{3}})^{\frac{3}{4}}+C(\int\limits_{T^{n-1}\times
S^{1}} |\Psi|^{4})^{\frac{1}{4}}. \eee
Here, $C$ also depends on $\Lambda$. Note that by assumption
\eqref{3.6} and the definition of $\epsilon_{1}$,
\be \omega:=\underset{t\in[T_{1}-1,T_{2}]}{\rm
sup}\int\limits_{[t,t+1]\times S^{1}}|d\Phi|^{2}+|\Psi|^{4} \leq
\epsilon_{1}\leq 1.  \nn \ee
Hence, by applying the $\epsilon$-regularity theorem to
\eqref{3.12}, we can conclude
\be \int\limits_{P^{n}}|\Psi|^{4} \leq C\omega^{2}+C\omega \leq
C\omega. \nn \ee
Similarly, we have
\be \int\limits_{P^{n}}|\nabla\Psi|^{\frac{4}{3}} \leq
C\omega^{\frac{1}{3}}. \nn \ee
Summing up the estimates on $P^{n}$ gives
\be \label{3.13} \int\limits_{P_{T_{1},T_{2}}}|\Psi|^{4}=
\sum\limits_{n=1}^{N_{0}}\int\limits_{P^{n}}|\Psi|^{4} \leq
CN_{0}\omega \leq C(\Lambda)\omega, \ee
and
\be \label{3.14}
\int\limits_{P_{T_{1},T_{2}}}|\nabla\Psi|^{\frac{4}{3}}=
\sum\limits_{n=1}^{N_{0}}
\int\limits_{P^{n}}|\nabla\Psi|^{\frac{4}{3}} \leq
CN_{0}\omega^{\frac{1}{3}} \leq C(\Lambda)\omega^{\frac{1}{3}}. \ee

In order to estimate $\int_{P_{T_{1},T_{2}}}|\Phi_{\theta}|^{2}$, we
again separate $P_{T_{1},T_{2}}$ into smaller parts as follows:
\be P_{T_{1},T_{2}}=\bigcup\limits_{n=1}^{N_{1}}\widetilde{P}^{n},
\widetilde{P}^{n}:=[\widetilde{T}^{n-1},\widetilde{T}^{n}]\times
S^{1}, \widetilde{T}^{0}=T_{1},
\widetilde{T}^{n}=\widetilde{T}^{n-1}+1,
\widetilde{T}^{N_{1}}=T_{2}. \nn \ee
Note that here, $N_{1}$ depends on $T_{1},T_{2}$. Since
$\widetilde{T}^{n}=\widetilde{T}^{n-1}+1$, by assumption \eqref{3.6}
and the $\epsilon$-regularity theorem, we have $|\Phi-q|\leq
C\omega^{\frac{1}{2}}$ on each part $\widetilde{P}^{n}$. Applying
Lemma \ref{lem3.2} on each part $\widetilde{P}^{n}$ and summing up
the inequalities gives
\bee \label{3.15} \int\limits_{P_{T_{1},T_{2}}}|\Phi_{\theta}|^{2}
&\leq& C\omega^{\frac{1}{2}}\int\limits_{P_{T_{1},T_{2}}}|d\Phi|^{2}
+C\omega^{\frac{1}{2}}\int\limits_{P_{T_{1},T_{2}}}|\Psi|^{4}
+\int\limits_{S_{1}}-\int\limits_{S_{2}}(\Phi-q)\Phi_{t}d\theta.
\eee
Since $q$ is equal to the mean value of $\Phi$ on $S_{i}$, by the
Poincar\'{e} inequality on $S_{i}$ and by H\"{o}lder's inequality,
it is easy to verify that
\bee \label{3.16} \int\limits_{S_{i}}|(\Phi-q)\cdot\Phi_{t}| \leq C
\int\limits_{S_{i}}|d\Phi|^{2}. \eee
By the $\epsilon$-regularity theorem and the Sobolev imbedding
theorem, we have
\bee \label{3.17} \int\limits_{S_{i}}|d\Phi|^{2} \leq C \omega .
\eee
Combining \eqref{3.15}, \eqref{3.16} and \eqref{3.17} gives
\bee \label{3.18} \int\limits_{P_{T_{1},T_{2}}}|\Phi_{\theta}|^{2}
\leq C\Lambda\omega^{\frac{1}{2}}+C\omega  \leq
C(\Lambda)\omega^{\frac{1}{3}}. \eee
Finally, by combining \eqref{3.18} with \eqref{3.13} and
\eqref{3.14}, we obtain
\be \int\limits_{P_{T_{1},T_{2}}}|\Phi_{\theta}|^{2}
+\int\limits_{P_{T_{1},T_{2}}}|\Psi|^{4}+\int\limits_{P_{T_{1},T_{2}}}|\nabla\Psi|^{\frac{4}{3}}
\leq C(\Lambda)\omega^{\frac{1}{3}}. \nn \ee
Thus we have finished the proof of Lemma \ref{lem3.3}.   \hfill
$\square$

\begin{lemma} \label{lem3.4} Let $(\Phi,\Psi)\in
C^{\infty}(P_{T_{1},T_{2}},N)$ be a Dirac-harmonic map. Then for
$t\in [T_{1},T_{2}]$,
\be \int_{\{t\}\times S^{1}}T(\Phi,\Psi)d\theta  \nn \ee
is independent of  $t\in [T_{1},T_{2}]$, where
\bee \label{3.19} T(\Phi,\Psi)=
(|\Phi_{t}|^2-|\Phi_{\theta}|^{2}-2i\phi_{t}\cdot
\Phi_{\theta})+({\rm Re} \langle\Psi,\partial_{t}\cdot\Psi_{t}
\rangle-i{\rm Re} \langle\Psi,\partial_{t}\cdot\Psi_{\theta}\rangle)
\eee
and $T(\Phi,\Psi)(dt+id\theta)^{2}$ is the holomorphic quadratic
differential of $(\Phi,\Psi)$ on $P_{T_{1},T_{2}}$.
\end{lemma}

\noindent\emph{Proof}. By Proposition \ref{pro2.2}, $T(\Phi,\Psi)$
is holomorphic on $P_{T_{1},T_{2}}$. The rest of the proof is
analogous to the case of harmonic maps (see Lemma 3.3 in \cite{Z2}).
\hfill $\square$

\begin{definition} \label{def3.1} Let $(\Phi,\Psi)\in
C^{\infty}(P_{T_{1},T_{2}},N)$ be a Dirac-harmonic map. Then we
define a complex number
\be \label{3.20} \alpha(\Phi,\Psi,P_{T_{1},T_{2}}):=
\int_{\{t\}\times S^{1}}T(\Phi,\Psi)d\theta \in \mathbb{C}  \ee
that is associated to $(\Phi,\Psi)$ along the cylinder
$P_{T_{1},T_{2}}$.
\end{definition}

\begin{remark} By Lemma \ref{lem3.4}, it is easy to verify that $\alpha(\Phi,\Psi,P_{T_{1},T_{2}})$ is
well-defined. Moreover, $\forall t_{1}<t'_{1}<t'_{2}<t_{2}$,
$\alpha(\Phi,\Psi, P_{t'_{1},t'_{2}})=\alpha(\Phi,\Psi,
P_{t_{1},t_{2}})$.
\end{remark}

\begin{lemma} \label{lem3.5} Let $(\Phi,\Psi)\in C^{\infty}(P_{T_{1},T_{2}},N)$
be a Dirac-harmonic map with

\noindent$\alpha$$=\alpha(\Phi,\Psi,P_{T_{1},T_{2}})$ and
\be
\label{3.21}\int\limits_{P_{T_{1}-1,T_{2}+1}}|d\Phi|^{2}+|\Psi|^{4}
\leq \Lambda < \infty. \ee
Then we have

\noindent(1)
\bee \label{3.22} |E(\Phi,P_{T_{1},T_{2}})-|{\rm
Re}\alpha|(T_{2}-T_{1})| &\leq&
2\int\limits_{P_{T_{1},T_{2}}}|\Phi_{\theta}|^{2}
+C\int\limits_{P_{T_{1},T_{2}}}|\Psi|^{4}
+C\int\limits_{P_{T_{1},T_{2}}}|\nabla\Psi|^{\frac{4}{3}} \nn\\
&&{}+C(\Lambda)(\int\limits_{P_{T_{1},T_{2}}}|\Psi|^{4})^{\frac{1}{2}}.\eee
\noindent(2)
\bee \label{3.23} |{\rm Im}\alpha|(T_{2}-T_{1})&\leq& C(\Lambda)
(\int\limits_{P_{T_{1},T_{2}}}|\Phi_{\theta}|^{2})^{\frac{1}{2}}
+C\int\limits_{P_{T_{1},T_{2}}}|\Psi|^{4}+C\int\limits_{P_{T_{1},T_{2}}}|\nabla\Psi|^{\frac{4}{3}}\nn\\
&&{}+C(\Lambda)(\int\limits_{P_{T_{1},T_{2}}}|\Psi|^{4})^{\frac{1}{2}}.\eee
Here $C, C(\Lambda)$ are constants independent of $T_{1},T_{2}$, and
$C(\Lambda)$ depends on $\Lambda$.
\end{lemma}

\noindent\emph{Proof}. By Definition \ref{def3.1} and \eqref{3.19},
we have
\bee \label{3.24} {\rm Re} \alpha &=&
\int_{0}^{2\pi}|\Phi_{t}|^{2}d\theta -
\int_{0}^{2\pi}|\Phi_{\theta}|^{2}d\theta - \int_{0}^{2\pi}{\rm
Re}\langle\Psi,\partial_{t}\cdot\Psi_{t}\rangle d\theta,                        \\
\label{3.25} {\rm Im}\alpha &=& -2\int_{0}^{2\pi}\Phi_{t}\cdot
\Phi_{\theta}d\theta - 2\int_{0}^{2\pi}{\rm Re}
\langle\Psi,\partial_{t}\cdot\Psi_{\theta}\rangle d\theta. \eee
Hence,
\bee \label{3.26} E(\Phi,P_{T_{1},T_{2}}) &=&
\int_{T_{1}}^{T_{2}}\int_{0}^{2\pi}|\Phi_{\theta}|^{2}d\theta dt
+\int_{T_{1}}^{T_{2}}\int_{0}^{2\pi}|\Phi_{t}|^{2}d\theta dt,   \nn  \\
&=& 2\int\limits_{P_{T_{1},T_{2}}}|\Phi_{\theta}|^{2}+ {\rm Re}
\alpha \cdot(T_{2}-T_{1})+\int\limits_{P_{T_{1},T_{2}}}{\rm
Re} \langle\Psi,\partial_{t}\cdot\Psi_{t}\rangle,  \nn  \\
&=& {\rm Re} \alpha
\cdot(T_{2}-T_{1})+\int\limits_{P_{T_{1},T_{2}}}(2|\Phi_{\theta}|^{2}+
{\rm Re} \langle\Psi,\partial_{t}\cdot\Psi_{t}\rangle). \eee
Let $a={\rm Re}\alpha \cdot(T_{2}-T_{1}), \quad
b=\int\limits_{P_{T_{1},T_{2}}}(2|\Phi_{\theta}|^{2}+{\rm Re}
\langle\Psi,\partial_{t}\cdot\Psi_{t}\rangle)$. Then by the
following inequality:
\bee |(a+b) - |a|| \leq |b|, \quad \forall a,b, a+b \geq 0 \nn \eee
we have
\bee  |E(\Phi,P_{T_{1},T_{2}})-|{\rm Re}\alpha|\cdot(T_{2}-T_{1})|
&\leq& \int\limits_{P_{T_{1},T_{2}}}(2|\Phi_{\theta}|^{2}+|{\rm Re}
\langle\Psi,\partial_{t}\cdot\Psi_{t}\rangle|)                    \nn  \\
&\leq& 2\int\limits_{P_{T_{1},T_{2}}}|\Phi_{\theta}|^{2}
+\int\limits_{P_{T_{1},T_{2}}}|\Psi|\cdot|\widetilde{\nabla}\Psi|\nn\\
&\leq&2\int\limits_{P_{T_{1},T_{2}}}|\Phi_{\theta}|^{2}
+C\int\limits_{P_{T_{1},T_{2}}}|\Psi|\cdot(|\nabla\Psi|+|d\Phi|\cdot|\Psi|)\nn\\
&\leq& 2\int\limits_{P_{T_{1},T_{2}}}|\Phi_{\theta}|^{2}
+C\int\limits_{P_{T_{1},T_{2}}}|\Psi|^{4}
+C\int\limits_{P_{T_{1},T_{2}}}|\nabla\Psi|^{\frac{4}{3}} \nn\\
&&{}+C(\Lambda)(\int\limits_{P_{T_{1},T_{2}}}|\Psi|^{4})^{\frac{1}{2}}.
\nn \eee
Here, in the last step, we used the Cauchy inequality, \eqref{3.21}
and the following inequality:
\bee ab\leq\frac{a^{4}+3b^{\frac{4}{3}}}{4}, \quad  \forall
a,b\geq0. \nn \eee
$C, C(\Lambda)$ are constants independent of $T_{1},T_{2}$, and
$C(\Lambda)$ depends on $\Lambda$.

By a similar argument, we can prove \eqref{3.23}.   \hfill $\square$

Now we consider a sequence of Dirac-harmonic maps from long spin
cylinders under certain assumptions. The following proposition gives
a ``bubble domain and neck domain" decomposition for such sequences,
which is analogous to the case of harmonic maps (c.f. Proposition
3.1 in \cite{Z2}).

\begin{proposition} \label{pro3.1} Let $(\Phi_{n},\Psi_{n})\in
C^{\infty}(P_{n},N)$ be a sequence of Dirac-harmonic maps with
$\alpha_{n}:=\alpha(\Phi_{n},\Psi_{n},P_{n})$, where
$P_{n}=[T_{n}^{1},T_{n}^{2}]\times S^{1}$ equipped with the
nontrivial spin structure. Assume that:
\begin{itemize}
\item[(1)] ``Long cylinder property"
\be  \label{3.27} 1 \ll T_{n}^{1} \ll T_{n}^{2},  \quad  i.e., \lim
_{n\rightarrow \infty} \frac{1}{T_{n}^{1}}= 0, \lim _{n\rightarrow
\infty} \frac{T_{n}^{1}}{T_{n}^{2}} = 0, \ee

\item[(2)] ``Uniform energy bound"
\be \label{3.28} E(\Phi_{n},\Psi_{n},P_{n}) \leq \Lambda < \infty,
\ee

\item[(3)] ``Asymptotic boundary conditions"
\be \label{3.29} \lim_{{n\rightarrow \infty}}
\omega(\Phi_{n},\Psi_{n},P_{T_{n}^{1},T_{n}^{1}+R})  =
\lim_{{n\rightarrow \infty}}
\omega(\Phi_{n},\Psi_{n},P_{T_{n}^{2}-R,T_{n}^{2}})
 = 0,  \forall R \geq 1, \ee
where
\be \omega(\Phi,\Psi,P_{T_{1},T_{2}})
:=\underset{t\in[T_{1},T_{2}-1]}{{\rm sup}}\int_{[t,t+1]\times
S^{1}}|d\Phi|^{2}+|\Psi|^{4}. \nn \ee
\end{itemize}
Then, after selection of a subsequence of
$(\Phi_{n},\Psi_{n},P_{n})$, either
\begin{itemize}
\item[(I)]
\be  \label{3.30} \lim_{{n\rightarrow \infty}}
\omega(\Phi_{n},\Psi_{n},P_{n})=0 \ee
or
\item[(II)]
$\exists K > 0$ independent of $n$ and $2K$ sequences
$\{a_{n}^{1}\}, \{b_{n}^{1}\}, \{a_{n}^{2}\}, \{b_{n}^{2}\},...,
\{a_{n}^{K}\}, \{b_{n}^{K}\}$ such that
\be \label{3.31} T_{n}^{1} \leq a_{n}^{1} \ll b_{n}^{1} \leq
a_{n}^{2} \ll
 b_{n}^{2} \leq ... \leq a_{n}^{K} \ll b_{n}^{K} \leq T_{n}^{2}, \ee
and
\be \label{3.32} (b_{n}^{i}-a_{n}^{i}) \ll T_{n}^{2}, \quad
i=1,2,...,K. \ee
Denote
\be J_{n}^{j}:=[a_{n}^{j},b_{n}^{j}]\times S^{1}, \quad j=1,2,...,K,
\nn \ee
\be I_{n}^{0}:= [T_{n}^{1}, a_{n}^{1}]\times S^{1}, I_{n}^{K}:= [
b_{n}^{K}, T_{n}^{2}]\times S^{1}, I_{n}^{i}:= [b_{n}^{i},
a_{n}^{i+1}]\times S^{1}, i=1,2,..., K-1. \nn \ee
Then
\begin{itemize}
\item[(i)] $\forall i=0,1,...K, \lim_{{n\rightarrow \infty}}
\omega(\Phi_{n},\Psi_{n},I_{n}^{i}) = 0.$ The maps
$(\Phi_{n},\Psi_{n}) :I_{n}^{i}\rightarrow N$ are necks
corresponding to collapsing homotopically nontrivial curves.

\item[(ii)] $\forall j=1,2,...,K,$ there are finitely many
Dirac-harmonic maps $(\omega^{j,l},\zeta^{j,l}): S^{2}\rightarrow N,
l=1,2,...,L_{j}$, such that:
\bee  \label{3.33} \lim \limits_{n \rightarrow
\infty}E(\Phi_{n},J_{n}^{j})&=& \sum_{l=1}^{L_{j}}E(\omega^{j,l}),  \\
\label{3.34} \lim \limits_{n \rightarrow
\infty}E(\Psi_{n},J_{n}^{j})&=& \sum_{l=1}^{L_{j}}E(\zeta^{j,l}).
\eee
\end{itemize}
\end{itemize}
\end{proposition}

\noindent\emph{Proof}. Analogous to the proof of Proposition 3.1 in
\cite{Z2}. One should be careful about the spin structures involved.
During the ``bubbling" procedure, we can obtain Dirac-harmonic maps
from $\mathbb{R}\times S^{1}$ which correspond to collapsing
homotopically nontrivial simple closed curves on $P_{n}$. By our
assumption that the spin structures on $P_{n}$ are nontrivial, the
induced spin structure on each domain $\mathbb{R}\times S^{1}$ is
nontrivial and thus can be conformally compactified to the unique
one on $S^{2}$. Hence, by using the conformal invariance and the
removability of singularities for Dirac-harmonic maps, we can get
Dirac-harmonic maps from $S^{2}$. For Dirac-harmonic maps from
$\mathbb{R}\times S^{1}$ corresponding to collapsing homotopically
trivial simple closed curves, see Remark \ref{rem2.3}. The energy
identities \eqref{3.33}, \eqref{3.34} follow from Theorem
\ref{thm2.3}. Note that here, for each $j$, we do not know whether
the maps $\omega^{j,l}, l=1,2,...,L_{j}$ are connected or not, since
the bubble tree convergence of Dirac-harmonic maps from a fixed spin
surface is still open. \hfill $\square$

The next lemma gives the asymptotics of the total energy of the
Dirac-harmonic necks $(\Phi_{n},\Psi_{n}): I_{n}^{i}\rightarrow N,
i=0,1,...,K$ as $n \rightarrow \infty$.

\begin{mainproposition} \label{mainpro3.1} Assumptions and
notations as in Proposition \ref{pro3.1}. W.l.o.g., we assume that
the limit $\lim \limits_{n \rightarrow \infty}|{\rm Re} \alpha
_{n}|\cdot|P_{n}|$ exists in $[0,+\infty]$, where
$|P_{n}|=T_{n}^{2}-T_{n}^{1}$. Then we have
\bee \label{3.35} \lim \limits_{n \rightarrow
\infty}\sum_{i=0}^{K}E(\Phi_{n},I_{n}^{i})&=& \lim \limits_{n
\rightarrow \infty}|{\rm Re} \alpha _{n}|\cdot|P_{n}|, \\
 \label{3.36}  \lim \limits_{n \rightarrow
\infty}\sum_{i=0}^{K}E(\Psi_{n},I_{n}^{i})&=&0. \eee
\end{mainproposition}

\noindent\emph{Proof}. We write
\bee  \label{3.37} \sum_{i=0}^{K}E(\Phi_{n},I_{n}^{i}) &=&
\sum_{i=0}^{K}|{\rm Re} \alpha _{n}|\cdot |I^{i}_{n}|
+ \sum_{i=0}^{K}(E(\Phi_{n},I_{n}^{i}) - |{\rm Re} \alpha _{n}|\cdot |I^{i}_{n}|)     \nn  \\
&=& I+II, \eee
where
\bee \label{3.38} I &:=& \sum_{i=0}^{K}|{\rm Re} \alpha _{n}|\cdot |I^{i}_{n}|                 \nn  \\
&=& |{\rm Re} \alpha _{n}|\cdot[(T_{n}^{2}-T_{n}^{1})-
\sum_{i=1}^{K}(b_{n}^{i}-a_{n}^{i})]                     \nn  \\
&=& |{\rm Re} \alpha _{n}|\cdot
(T_{n}^{2}-T_{n}^{1})\cdot(\frac{T_{n}^{2}}{T_{n}^{2}-T_{n}^{1}})\cdot
[(1-\frac{T_{n}^{1}}{T_{n}^{2}})-
\sum_{i=1}^{K}\frac{(b_{n}^{i}-a_{n}^{i})}{T_{n}^{2}}] \eee
and
\bee II := \sum_{i=0}^{K}(E(\Phi_{n},I_{n}^{i}) - |{\rm Re} \alpha
_{n}|\cdot |I^{i}_{n}|), \nn \eee
By Lemma \ref{lem3.3}, Lemma \ref{lem3.5} and Proposition
\ref{pro3.1}, for $n$ large enough, we have
\bee |II|&\leq& \sum_{i=0}^{K}|E(\Phi_{n},I_{n}^{i}) - |{\rm Re}
\alpha
_{n}|\cdot |I^{i}_{n}||      \nn  \\
&\leq& \sum_{i=0}^{K}
(2\int\limits_{I_{n}^{i}}|(\Phi_{n})_{\theta}|^{2}
+C\int\limits_{I_{n}^{i}}|\Psi_{n}|^{4}
+C\int\limits_{I_{n}^{i}}|\nabla\Psi_{n}|^{\frac{4}{3}}
+C(\Lambda)(\int\limits_{I_{n}^{i}}|\Psi_{n}|^{4})^{\frac{1}{2}})
\nn \eee
Note that in Proposition \ref{pro3.1}, after passing to
subsequences, the local energy of $(\Phi_{n},\Psi_{n})$ over a small
neighborhood of the two boundary components of $I^{i}_{n}$ can be
arbitrary small. Thus, by Lemma \ref{lem3.5} and Proposition
\ref{pro3.1},
\bee \label{3.39} |II| &\leq& C(\Lambda)\sum_{i=0}^{K}\left[
(\omega(\Phi_{n},\Psi_{n},I_{n}^{i}))^{\frac{1}{3}}+(\omega(\Phi_{n},\Psi_{n},I_{n}^{i}))^{\frac{1}{6}}\right]
\rightarrow 0, n\rightarrow \infty. \eee
From Proposition \ref{pro3.1}, we have $ 1\ll T_{n}^{2}\ll
T_{n}^{2}, 1\ll(b_{n}^{i}-a_{n}^{i}) \ll T_{n}^{2}, i=1,2,...,K.$
Hence, by \eqref{3.38}, we get $\lim \limits_{n \rightarrow \infty}
I=\lim \limits_{n \rightarrow \infty}|{\rm Re} \alpha
_{n}|\cdot|P_{n}|$. Now \eqref{3.35} follows immediately from
\eqref{3.37} and \eqref{3.39}.

Recall that $\lim_{n\rightarrow
\infty}\omega(\Phi_{n},\Psi_{n},I_{n}^{i})=0, i=0,1,...,K$. Thus,
for $n$ sufficiently large, applying Lemma \ref{lem3.3} on each
$I_{n}^{i}$, summing up the inequalities and finally taking the
limit ($n \rightarrow \infty$), we can prove \eqref{3.36}.  \hfill
$\square$

By similar arguments as in the proof of Main Proposition
\ref{mainpro3.1}, we get
\begin{proposition} \label{pro3.2} With the same assumptions and notations as Proposition \ref{pro3.1}, we have
\bee  \limsup \limits_{n \rightarrow \infty}|{\rm Re} \alpha
_{n}|\cdot|P_{n}|\leq \Lambda, \quad
 \limsup \limits_{n \rightarrow
\infty}|{\rm Im} \alpha _{n}|\cdot|P_{n}| = 0. \eee
\end{proposition}

\noindent\emph{Proof}. By Lemma \ref{lem3.3}, Lemma \ref{lem3.5},
Proposition \ref{pro3.1} and Main Proposition \ref{mainpro3.1}.
\hfill $\square$

\vskip 0.2cm

Combining the results obtained before, we state the following
generalized energy identities for Dirac-harmonic maps from long spin
cylinders:

\begin{theorem} \label{thm3.1} Assumptions and
notations as in Main Proposition \ref{mainpro3.1}. Then there are
finitely many Dirac-harmonic maps
$(\omega^{j,l},\zeta^{j,l}):S^{2}\rightarrow N,$ $l=1,2,...,L_{j}$;
$j=1,2,...,K$, such that after selection of a subsequence of
$(\Phi_{n},\Psi_{n},P_{n})$, the following holds:
\bee \label{3.41} \lim \limits_{n \rightarrow
\infty}E(\Phi_{n},P_{n})&=&\sum_{j=1}^{K}\sum_{l=1}^{L_{j}}E(\omega^{j,l})+
\lim \limits_{n \rightarrow \infty}|{\rm Re} \alpha _{n}|\cdot|P_{n}|,     \\
 \label{3.42} \lim\limits_{n\rightarrow\infty}E(\Psi_{n},P_{n})&=&
\sum\limits_{j=1}^{K}\sum\limits_{l=1}^{L_{j}}E(\zeta^{j,l}). \eee
\end{theorem}

\noindent\emph{Proof}. By Proposition \ref{pro3.1} and Main
Proposition \ref{mainpro3.1}.   \hfill $\square$

\vskip 0.2cm

As a corollary of Theorem \ref{thm3.1}, we have

\begin{corollary} \label{cor3.1} Assumptions and notations as in Proposition \ref{pro3.1}.
Then $(\Phi_{n},\Psi_{n},P_{n})$ subconverges in $W^{1,2}\times
L^{4}$ modulo bubbles, i.e., in the limit, the necks contain no
energy, if and only if
\bee \liminf\limits _{n\rightarrow \infty}|{\rm Re} \alpha
_{n}|\cdot |P_{n}|=0. \eee
\end{corollary}

\vskip 0.5cm

\section{Dirac-harmonic maps from degenerating spin surfaces}
\label{section 4}

In this section we will apply the results developed in the previous
sections to prove our main theorems stated in the introduction.

To begin the proofs, we shall give a brief introduction to
deformations of spin surfaces. A compact connected Riemann surface
with a spin structure can be viewed as an algebraic curve with a
theta characteristic, i.e., a square root of the canonical bundle
\cite{Mu}, \cite{A}. The moduli space of curves with theta
characteristics can be compactified algebraically by generalizing
the notion of theta characteristics to the case of singular curves
(a good reference is \cite{C}).

Here, following the discussions in \cite{B}, we present a geometric
topological description of the degeneration of spin surfaces. Let
$(M_{n},h_{n},c_{n},\mathfrak{S}_{n})$ be a sequence of closed
hyperbolic Riemann surfaces of genus $g>1$ with spin structures
$\mathfrak{S}_{n}$. We assume that $(M_{n},h_{n},c_{n})$ degenerates
to a hyperbolic Riemann surface $(M,h,c)$ by collapsing $p$ ($1\leq
p\leq 3g-3$) pairwise disjoint simple closed geodesics
$\{\gamma_{n}^{j},j=1,2,...,p\}$ on $M_{n}$. For each $j$, the
geodesics $\gamma_{n}^{j}$ degenerate into a pair of punctures
$(\mathcal{E}^{j,1},\mathcal{E}^{j,2})$. Let
$(\overline{M},\overline{c})$ be the normalization of $(M,h,c)$. Let
$\tau_{n}: M \rightarrow M_{n}\setminus\cup_{j=1}^{p}\gamma_{n}^{j}$
be the corresponding diffeomorphisms realizing the degeneration
(c.f. \cite{Z2}). For each $n$, the diffeomorphism $\tau_{n}$ and
the spin structure $\mathfrak{S}_{n}$ together determine a pull-back
spin structure on $M$. If we identify spin structures on $M$ with
elements of $H^{1}(M,\mathbb{Z}/2\mathbb{Z})$, then a spin structure
on $M$ corresponds uniquely to a spin structure on $\overline{M}$
together with a choice of an even number of punctures along which
the spin structure is trivial; the induced spin structures along the
remaining punctures are nontrivial (c.f. Sect. 8 in \cite{B}). It is
clear that there are finitely many spin structures on a surface with
punctures. Thus, by taking subsequences, we can assume that
$\tau_{n}$ is compatible with the spin structures
$\mathfrak{S}_{n}$, namely the pull-back spin structure on the limit
surface $M$ is fixed. Let us denote it by $\mathfrak{S}$. In
particular, for each $j$, $\mathfrak{S}$ is nontrivial or trivial
along the pair of punctures $(\mathcal{E}^{j,1},\mathcal{E}^{j,2})$
if and only if $\mathfrak{S}_{n}$ is nontrivial or trivial along the
geodesic $\gamma_{n}^{j}$ for all $n$.

Now, we consider a sequence of smooth Dirac-harmonic maps
\bee (\phi_{n},\psi_{n}): (M_{n},h_{n},c_{n},\mathfrak{S}_{n})
\rightarrow N  \nn \eee
with uniformly bounded energy $E(\phi_{n},\psi_{n},M_{n})\leq
\Lambda <\infty$, where $(M_{n},h_{n},c_{n},\mathfrak{S}_{n})$ is a
sequence of closed hyperbolic Riemann surfaces of genus $g>1$ with
spin structures $\mathfrak{S}_{n}$. We assume that
$(M_{n},h_{n},c_{n})$ degenerates to a hyperbolic Riemann surface
$(M,h,c)$ by collapsing $p$ ($1\leq p\leq 3g-3$) pairwise disjoint
simple closed geodesics $\{\gamma_{n}^{j},j=1,2,...,p\}$ on $M_{n}$.
Let $\alpha_{n}^{j}:=\alpha(\phi_{n},\psi_{n},P_{n}^{j})$ be the
quantities associated to $(\phi_{n},\psi_{n})$ along the $j$-th
cylindrical collar $P_{n}^{j}$ as in Definition \ref{def3.1}. By
taking subsequences, we can assume that the limits $$\lim\limits_{n
\rightarrow \infty}|{\rm Re} \alpha _{n}^{j}|\cdot
\frac{2\pi^{2}}{l_{n}^{j}}, \quad  j=1,2,...,p$$ exist in
$[0,\infty]$. After passing to a further subsequence, we assume that
the pull back of $\mathfrak{S}_{n}$ via $\tau_{n}$ is a fixed spin
structure $\mathfrak{S}$ on $M$. Note that $M$ has $p$ pairs of
punctures.

We require the following additional assumption:
\bee \label{4.1} \emph{All punctures of the limit spin surface
$(M,\mathfrak{S})$ are of Neveu-Schwarz type.} \eee
It is equivalent to say that the spin structure $\mathfrak{S}$ is
nontrivial along all punctures of $M$ and $ \mathfrak{S}_{n}$ is
nontrivial along all degenerating collars $P_{n}^{j}, j=1,2,...,p$.
Thus, the spin structure $\mathfrak{S}$ on $M$ extends to some spin
structure $\overline{\mathfrak{S}}$ on $\overline{M}$ as explained
before.

Replacing the data on $M_{n}$ by the pull-back data on $M$ and
passing to subsequences, we can fix the spinor bundle $\Sigma M$ and
think of the hyperbolic metrics and the compatible complex
structures $(h_{n},c_{n})$ as all living on the limit surface $M$
and converging in $C_{loc}^{\infty}$ to $(h,c)$. Let $\nabla_{n}$ be
the connection on $\Sigma M$ coming from $h_{n}$ and $\nabla$ the
connection on $\Sigma M$ coming from $h$. Then, we can consider
$(\phi_{n}, \psi_{n})$ as a sequence of Dirac-harmonic maps defined
on $(M,h_{n},c_{n},\mathfrak{S})$ with respect to
$(c_{n},\nabla_{n})$.

\vskip 0.2cm

\noindent\emph{Proof of Theorem \ref{thm1.2}}. Analogous to the
proof of Theorem 1.1 in \cite{Z2}. Note that here, $(h_{n},c_{n})$
converges in $C_{loc}^{\infty}$ to $(h,c)$, as $n\rightarrow
\infty$. Hence, the tensor $\nabla_{n}-\nabla$ converges to zero in
$C_{loc}^{\infty}$ away from the punctures of $M$ and the energy
functional corresponding to $\nabla_{n}$ over any compact subset of
$M$ is uniformly equivalent to the one corresponding to $\nabla$
over the same domain. By the uniform energy bound, we can assume
that $(\phi_{n},\psi_{n})$ subconverges weakly to a limit
$(\phi,\psi)$ in $W_{loc}^{1,2}(M,N)\times L_{loc}^{4}(\Sigma
M\otimes \mathbb{R}^{K})$ with respect to $(c,\nabla)$. By similar
arguments as in the proof of Theorem \ref{thm1.1}, there exist
finitely many blow-up points $\{x_{1},x_{2},...,x_{I}\}$ away from
the punctures, such that $(\phi_{n},\psi_{n})$ subconverges to
$(\phi,\psi)$ in $C_{loc}^{\infty}\times C_{loc}^{\infty}$ on
$M\setminus\{x_{1},x_{2},...,x_{I}\}$ with respect to $(c,\nabla)$.
Furthermore, $(\phi,\psi)$ is actually a smooth Dirac-harmonic map
defined on $M$ with respect to $(c,\nabla)$. Note that the complex
structure $c$ on $M$ extends to some complex structure
$\overline{c}$ on $\overline{M}$. By our assumption \eqref{4.1}, the
spin structure $\mathfrak{S}$ on $M$ extends to some spin structure
$\overline{\mathfrak{S}}$ on $\overline{M}$. Hence, applying the
removable singularity theorem, we have that $(\phi,\psi)$ extends to
a smooth Dirac-harmonic map defined on
$(\overline{M},\overline{c},\overline{\mathfrak{S}})$.

The energy concentration at the blow-up points
$\{x_{1},x_{2},...,x_{I}\}$ that are away from the punctures is
analogous to the case in Theorem \ref{thm1.1}. With similar
arguments as in the proof of Theorem \ref{thm1.1} in \cite{Z2}, the
energy concentration near the punctures can be reduced to the study
of Dirac-harmonic maps from degenerating collars $P_{n}^{j},
j=1,2,...,p$. By assumption \eqref{4.1}, for all $n$, the spin
structure $ \mathfrak{S}_{n}$ is nontrivial along each of the
degenerating collars $P_{n}^{j}, j=1,2,...,p$. Thus, we can apply
the results in Section 3, especially Theorem \ref{thm3.1}, to
capture the energy loss along the collars.    \hfill $\square$

\vskip 0.2cm

\noindent\emph{Proof of Corollary \ref{cor1.1}}. By Theorem
\ref{thm1.2} and Corollary \ref{cor3.1}.    \hfill $\square$

\vskip 0.2cm

\noindent\emph{Proof of Proposition \ref{pro1.1}}. By Theorem
\ref{thm1.2} and Proposition \ref{pro3.2}.       \hfill $\square$

\vskip 0.2cm

\begin{remark} The assumption \eqref{4.1} can be satisfied by
choosing suitable topological types of degeneration. For example, we
consider a closed hyperbolic Riemann surface $M$ of genus $g>1$. Let
$\gamma$ be a simple closed geodesic on $M$. If $M \setminus \gamma$
is disconnected, then any spin structure on $M$ must be nontrivial
along the collar around $\gamma$ (c.f. Sect. 8 in \cite{B}). In this
case, pinching $\gamma$ to a point and deleting it gives a pair of
Neveu-Schwarz punctures. If $M \setminus \gamma$ is connected, then
the collar around $\gamma$ can carry two spin structures \cite{B}.
Hence, both types of singularities can occur when pinching $\gamma$
to a point. More generally, if $M$ carries a spin structure
$\mathfrak{S}$ with Arf invariant $1$, then there is a collection of
pairwise disjoint simple closed curves $\gamma_{j}: S^{1}\rightarrow
M, j=1,2,...,g$ representing linearly independent elements in
$H_{1}(M,\mathbb{Z})$, such that $\mathfrak{S}$ is nontrivial along
each of the $\gamma^{j}$ (c.f. Sect. 3 in \cite{AB}).
\end{remark}

\vskip 0.5cm

\noindent Miaomiao Zhu\\
Max Planck Institute for Mathematics in the Sciences\\
Inselstr.\ 22-26, D-04103 Leipzig, Germany\\
E-mail: Miaomiao.Zhu@mis.mpg.de

\end{document}